\documentclass[11pt]{article}
\usepackage{times}
\usepackage{amssymb,amsmath}
\usepackage{epsf}
\usepackage{times,bm,mathrsfs}  
\usepackage{amsmath}
\usepackage{amssymb}
\usepackage{amsfonts}
\usepackage{mathrsfs}
\usepackage{bbm}
\usepackage{color}
\usepackage{graphicx}
\usepackage{amscd}
\usepackage{epsfig}

\definecolor{Red}{rgb}{1,0,0}
\definecolor{Blue}{rgb}{0,0,1}
\definecolor{Olive}{rgb}{0.41,0.55,0.13}
\definecolor{Green}{rgb}{0,1,0}
\definecolor{MGreen}{rgb}{0,0.8,0}
\definecolor{DGreen}{rgb}{0,0.55,0}
\definecolor{Yellow}{rgb}{1,1,0}
\definecolor{Cyan}{rgb}{0,1,1}
\definecolor{Magenta}{rgb}{1,0,1}
\definecolor{Orange}{rgb}{1,.5,0}
\definecolor{Violet}{rgb}{.5,0,.5}
\definecolor{Purple}{rgb}{.75,0,.25}
\definecolor{Brown}{rgb}{.75,.5,.25}
\definecolor{Grey}{rgb}{.5,.5,.5}
\definecolor{Black}{rgb}{0,0,0}

\newcommand{\br}{{\mbox{{\rm br}}}}

\def\path{{\tt path}}
\def\lam{{\lambda}}
\def\one{{\bf 1}}

\newcommand{\acal}{\mathcal{A}}
\newcommand{\bcal}{\mathcal{B}}
\newcommand{\ccal}{\mathcal{C}}

\newcommand{\tcal}{\mathcal{T}}

\newcommand{\A}{\ccal}

\newcommand{\eps}{\varepsilon}

\newcommand{\bdm}{\begin{displaymath}}
\newcommand{\edm}{\end{displaymath}}
\newcommand{\bea}{\begin{eqnarray*}}
\newcommand{\eea}{\end{eqnarray*}}
\newcommand{\bean}{\begin{eqnarray}}
\newcommand{\eean}{\end{eqnarray}}

\newcommand{\prob}{\mathbb{P}}
\newcommand{\expec}{\mathbb{E}}

\newcommand{\diff}{\mathrm{d}}

\renewcommand{\deg}{{ b}}
\newcommand{\cutsets}{\mathrm{cutsets}}
\newcommand{\pimp}{\pi_{-/+}}
\newcommand{\deltat}{\Delta}
\newcommand{\that}{\widehat{T}}

\newcommand{\yhat}{\hat y}
\newcommand{\yyhat}{\widehat{Y}}
\newcommand{\mom}[1]{{\bar #1}}
\newcommand{\momp}[1]{{{\bar #1}}_+}
\newcommand{\momm}[1]{{{\bar #1}_-}}

\newtheorem{theorem}{Theorem}
\newtheorem{proposition}{Proposition}

\newtheorem{definition}{Definition}

\newtheorem{lemma}{Lemma}

\newenvironment{proof}{\noindent{\textbf{Proof:}}}{$\blacksquare$\vskip\belowdisplayskip}

\author{
DRAFT - DO NOT CIRCULATE
}
\title{
The Kesten-Stigum Reconstruction Bound\\
Is Tight for Roughly Symmetric Binary Channels\\
}

\author{Christian Borgs\thanks{{ Microsoft Research.}} \and Jennifer Chayes\thanks{{ Microsoft Research.}} \and
Elchanan Mossel\thanks{Dept. of Statistics, U.C. Berkeley. Supported by an Alfred Sloan fellowship
  in Matheamtics and by NSF grants DMS-0528488, DMS-0504245 and
  DMS-0548249 (CAREER). Most of this work was done while
  visiting { Microsoft} Research.} \and
  Sebastien Roch\thanks{{ Dept. of Statistics, U.C. Berkeley.  Some 
  of this work was done while
  visiting Microsoft Research.}}}
\begin{document}

\maketitle

\thispagestyle{empty}

\begin{abstract}
{
We establish the exact threshold for the reconstruction problem for a
binary asymmetric channel on the $\deg$-ary tree, provided that the asymmetry
is sufficiently small.  This is the first exact reconstruction threshold
obtained in roughly a decade.  We discuss the implications of our result for
Glauber dynamics, phylogenetic reconstruction, and so-called ``replica symmetry
breaking'' in spin glasses and random satisfiability problems. }

\end{abstract}

\bigskip


\clearpage

\section{Introduction}
Let
\begin{equation} \label{eq:channel}
M
=
\frac{1}{2}\left[
\left(\begin{matrix}
1+\theta & 1-\theta\\
1-\theta & 1+\theta
\end{matrix}\right)
+
\delta\left(\begin{matrix}
-1 & 1\\
-1 & 1
\end{matrix}\right)
\right],
\end{equation}
be a binary asymmetric channel with second eigenvalue
$\theta$ and $T_{\deg}$ be a complete $\deg$-ary tree.
{
The ``reconstruction problem'' is the problem of determining
the state of the root, given the distribution of the
Markov chain on level $n$ of the tree, as $n$ gets larger and
larger (a precise definition is given below).  For the
symmetric binary channel ($\delta=0$), it was known since 1995 
that the reconstruction problem is solvable if and only
if $\deg\theta^2>1$. 
For all other channels, it was also known and easy to prove that
$\deg\theta^2>1$ implies solvability, but {\em exact
non-solvability} results were not known.
Here we show that
this bound is tight provided that $M$ is close enough to { symmetric---i.e.,} we
show that the reconstruction problem for $M$ on $T$ is not
solvable if $\deg\theta^2 \leq 1$ and $|\delta|$ is sufficiently
small.}

The reconstruction problem is intimately related to mixing of Glauber
dynamics and to { phylogenetic} reconstruction. Moreover, it was recently
claimed that { the reconstruction problem corresponds to the
``replica symmetry broken'' solution of the spin glass on the tree.
replica symmetry breaking is a central notion in the statistical
physics theory of spin glasses and random satisfiability problems.
We discuss potential applications of our results in these different areas.}

\subsection{Definitions and Main Result}
Let $T=(V,E,\rho)$ be a tree
$T$ with nodes $V$, edges $E$
and root $\rho \in V$. We direct all edges away from the root, so that
if $e=(x,y)$ then $x$ is on the
path connecting $\rho$ to $y$.
Let $d(\cdot,\cdot)$ denote the graph-metric distance on $T$, and
$L_n = \{v \in V: d(\rho,v) = n\}$ be the { $n^{\rm th}$} level of the tree.
For $x \in V$ and $e=(y,z) \in E$, we denote $|x| = d(\rho,x)$,
$d(x,(y,z)) = \max\{d(x,y),d(x,z)\}$, and
$|e| = d(\rho,e)$.
The {\em $\deg$-ary}  tree is the infinite rooted tree where each vertex
has exactly  $\deg$ children.

A Markov chain on the tree $T$ is a probability measure { defined on
the state space}
$\A^V$, where $\A$ is a finite set. 
{ Assume first that $T$ is finite and, for each edge $e$ of $T$, let $M^e = (M^e_{i,j})_{i,j \in \A}$
be a stochastic matrix.} In this case the probability measure defined
by $(M^e : e \in E)$ on $T$ is given by
\begin{equation} \label{eq:finite_def}
\overline{\mu}_{\ell}(\sigma) =
\one_{\{\sigma({\rho}) = \ell\}} \prod_{(x,y) \in E}
M^{(x,y)}_{\sigma(x),\sigma(y)}.
\end{equation}
In other words, the root state $\sigma({\rho})$ satisfies
$\sigma({\rho}) = \ell$ and then
each vertex iteratively chooses its state from the one of its parent by
an application of the Markov transition rule given by $M^e$ (and all
such applications
are independent).
We can define the measure $\overline{\mu}_{\ell}$
on an infinite tree as well, by Kolmogorov's extension theorem,
but we will not need chains on infinite trees in this paper (see
\cite{Georgii:88} for basic properties of Markov chains on trees).

Instead, for an infinite tree $T$, we let $T_n = (V_n,E_n,\rho)$, where
$V_n = \{x \in V: d(x,\rho) \leq n\},
  E_n = \{e \in E:  d(e,\rho) \leq n\}$ and define
  $\overline{\mu}_{\ell}^n$ by (\ref{eq:finite_def})  for $T_n$.
We are particularly interested in the distribution of the states
$\sigma(x)$ for $x\in L_n$, the set of leaves in $T_n$.
This distribution,
denoted by $\mu_k^n$,
is the projection of
$\overline{\mu}_k^n$  on $\A^{L_n}$
given by
\begin{equation} \label{eq:proj}
\mu_k^n(\sigma) =
\sum_{\bar\sigma} \{\overline{\mu}_k^n({\overline{\sigma}}) :
\overline{\sigma} | L_n = \sigma\}.
\end{equation}

Recall that for distributions $\mu$ and $\nu$ on the same space
$\Omega$ the total variation distance between $\mu$ and $\nu$ is
\begin{equation} \label{eq:totalvar}
D_V(\mu,\nu) = \frac{1}{2} \sum_{\sigma \in \Omega}
|\mu(\sigma) - \nu(\sigma)|.
\end{equation}

\begin{definition}[Reconstructibility]\label{def:reconstruction}
The reconstruction problem for the infinite tree $\tcal$ and $(M^e : e \in E)$ is
{\bf solvable} if there exist $i,j \in \A$ for which
\begin{equation} \label{cond:arb_l1}
\liminf_{n \to \infty} D_V(\mu_i^n,\mu_j^n) > 0.
\end{equation}
When $M^e = M$ for all $e$, we say that the reconstruction problem is
solvable for $T$ and $M$.
\end{definition}
We will be mostly interested in binary channels,
i.e., transition matrices on the state space
$\{\pm \}$. In this case, the definition above says that the reconstruction
problem is solvable if
\begin{equation} \label{cond:arb_l12}
\liminf_{n \to \infty} D_V(\mu_{+}^n,\mu_{-}^n) > 0.
\end{equation}

{ Our main result is the following:}

\begin{theorem}[Main Result] \label{thm:d1}
{ For all $\deg \geq 2$, there exists a $\delta_0 > 0$
such that for all $|\delta|\leq \delta_0$,
the reconstruction problem for $M$ on the $\deg$-ary tree $T_{\deg}$ is not solvable if
$\deg \theta^2 \leq 1$.}
\end{theorem}

\subsection{Previous Results}

The study of the reconstruction problem began in the
seventies~\cite{Spitzer:75,Higuchi:77} when the problem was
introduced in terms of the extremality of the free Gibbs measure on
the tree. In~\cite{Higuchi:77} it is shown that the reconstruction
problem for the binary symmetric channel (equation~(\ref{eq:channel}) where
$\delta = 0$) on the binary tree is solvable when $2 \theta^2 > 1$.
{ This in fact} follows from a previous work~\cite{KestenStigum:66}
which implies that for any Markov chain $M$, the reconstruction
problem on the $\deg$-ary tree is solvable if $\deg \theta^2 > 1$ where
$\theta$ is { the second largest eigenvalue of
$M$ in absolute value.}

Proving { non-reconstructibility} turned out to be harder. While
coupling { arguments} easily yield non-reconstruction, these
arguments are typically not tight.
{ A natural way to try to prove non-reconstructibility is to analyze
recursions 1) in terms of random variables each of whose values
is the expectation
of the chain at a vertex, given the state at the leaves of the subtree below it, 2) in
terms of ratios of such probabilities, or 3) in terms of log-likelihood
ratios of such probabilities.}
Such recursions were analyzed for a closely related model in~\cite{CCST:86}.
Both the reconstruction model and the model
analyzed in~\cite{CCST:86} deal with the correlation between 
the $n^{\mathrm{th}}$-level and
the root. However, while in the reconstruction problem, the two
random variables are generated according to the Markov model on the
tree, in~\cite{CCST:86} the nodes at level $n$ are set to have an
i.i.d. distribution and the root has the conditional distribution thus
induced.

In spite of this important difference, the two models are closely
related. In particular, in~\cite{CCST:86} it is shown that for the
binary tree, the
correlation between level $n$ and the root decays if and only if $2
\theta^2 \leq 1$.
Building on the techniques of~\cite{CCST:86} it was finally shown
in~\cite{BlRuZa:95} that the reconstruction problem for the binary
symmetric channel is solvable if and only if $2 \theta^2 > 1$.
This result was later reproven in various
ways~\cite{EvKePeSc:00,Ioffe:96a,BeKeMoPe:05,MaSiWe:04}.

The elegance of the threshold $\deg \theta^2 = 1$ raised the hope that
it is the threshold for reconstruction for general channels.
However, { previous attempts to generalize any} of the proofs to
other channels have failed. Moreover in~\cite{Mossel:01} it was shown that for
asymmetric binary channels and for symmetric channels on large
alphabets the reconstruction problem is solvable in cases where
$\deg \theta^2 < 1$. In fact~\cite{Mossel:01} contains an example of a
channel satisfying $\theta = 0$ for which the reconstruction problem
is solvable.
On the other hand, in~\cite{MosselPeres:03,JansonMossel:04} it is
shown that
the threshold $\deg \theta^2 = 1$ is the threshold for two variants of
the reconstruction problem, ``census reconstruction'' and ``robust
reconstruction''.

The results above led some to believe that ``reconstruction''
unlike its siblings ``census reconstruction'' and ``robust
reconstruction'' is an extremely sensitive property and that the
threshold $\deg \theta^2 = 1$ is tight only for { the} binary symmetric
channel.
This conceptual picture was shaken by recent results in { the theoretical
physics literature}~\cite{MezardMontanari:06} where using variational principles
developed in the context of ``replica symmetry breaking'' it is
suggested that the bound $\deg \theta^2 = 1$ is tight for symmetric
channels on $3$ and (maybe) $4$
letters.

In Theorem~\ref{thm:d1} we give the first tight threshold for the
reconstruction problem for channels other than binary symmetric
channels. We show that for asymmetric { channels} that are close to
symmetric, the Kesten-Stigum bound $\deg \theta^2 = 1$ is tight for reconstruction.
Our proof builds on ideas
from~\cite{CCST:86,BlRuZa:95,EvKePeSc:00,PemantlePeres:06} and is
extremely simple. { In addition to giving a new result for the asymmetric
channel, our proof also provides a much simpler proof of the previously
known result for the binary {\it symmetric} channel.}

\subsection{The Reconstruction Problem in Mixing, Phylogeny and
  Replicas}

\paragraph{Mixing of Markov Chains.}
One of the main themes at the intersection of statistical physics
and theoretical computer science
in recent years has been the study of connections between spatial and
temporal mixing. It is widely accepted that spatial mixing and
temporal mixing of dynamics go hand in hand though this was proven
only in restricted settings.

In particular, the spatial mixing condition is usually stated in terms
of uniqueness of Gibbs measures. However, as shown
in~\cite{BeKeMoPe:05} , this spatial condition is too strong. In
particular, it is shown in~\cite{BeKeMoPe:05} that the spectral gap of
continuous-time Glauber dynamics for the Ising model with no external
field and no boundary conditions
on the $\deg$-ary tree is $\Omega(1)$ whenever $\deg \theta^2 <
1$. This should be compared with the uniqueness condition on the tree
given by $\deg \theta < 1$. In~\cite{MaSiWe:04} this result is extended
to the log Sobolev constant. In~\cite{MaSiWe:04} it is also shown that
for measures on trees, a super-linear decay of point-to-set
correlations implies an $\Omega(1)$ spectral gap for the Glauber
dynamics with free boundary conditions.

Thus our { results  not} only give the exact threshold for
reconstructibility. They also yield an exact threshold for mixing of
Glauber dynamics on the tree for Ising models with a small external
field. The details are omitted from this extended abstract.

\paragraph{Phylogenetic Reconstruction.}
Phylogenetic reconstruction is a major task of systematic
biology~\cite{Felsenstein:04}.
It was recently shown in~\cite{DaMoRo:06} that for binary symmetric
channels, also called CFN models in evolutionary biology,
 the sampling efficiency of phylogenetic reconstruction is determined
by the reconstruction threshold.
Thus if for all edges of the tree it holds that $2 \theta^2 > 1$ the
tree can be recovered efficiently from $O(\log n)$ samples. If $2
\theta^2 < 1$, then~\cite{Mossel:04} implies that $n^{\Omega(1)}$
samples are needed.
In fact, the proof of the lower bound in~\cite{Mossel:04} implies the
lower bound $n^{\Omega(1)}$ whenever the reconstruction problem is
{\em exponentially} unsolvable. In other words, if
$\liminf_{r \to \infty} D_V(\mu_{+}^r,\mu_{-}^r) = \exp(-\Omega(r))$
then a lower bound of $n^{\Omega(1)}$ holds for { phylogenetic}
reconstruction.

Thus, our results here imply $n^{\Omega(1)}$ lower bounds
for { phylogenetic}
reconstruction for asymmetric channels such that $2 \theta^2 < 1$ and
$|\delta| < \delta_0$. The details are omitted from this extended abstract.
It is natural to conjecture that this is tight and that if $2 \theta^2
> 1$ then phylogenetic reconstruction may be achieved with $O(\log n)$
sequences.

\paragraph{Replica Symmetry Breaking.}
The replica and cavity methods were invented in the theoretical physics
literature to solve  Ising spin glass problems on the complete graph---the
so-called Sherrington-Kirkpatrick model. These
methods, while not mathematically rigorous, led to
numerous predictions on the spin glass and other models on
dense graphs, a few of which were proved many years later.  When applied to
random satisfiability problems, which turn out to be equivalent
to dilute spin glasses---i.e., spin glasses on sparse random
 graphs---these methods  led to the
empirically best algorithms for solving random satisfiability
problems~\cite{MezardZecchina:02,MePaZe:03}.

A central concept in this theory is the notion of a ``glassy phase''
of { the spin glass measure. In the glassy phase,
the distribution on the
random graphs decomposes into an exponential number of ``lumps''.
One of the standard techniques for determining the glassy phase is via
``replica symmetry breaking''.  Moreover, there
are certain glassy phases for which the
replica symmetry breaking is relatively simple---those which are said to
have ``one-step replica symmetry breaking''; and others in which the replica
symmetry breaking is more complicated---those with so-called ``full replica
symmetry breaking''.}

In a recent paper~\cite{MezardMontanari:06} it is claimed that the
parameters for which a ``glassy phase occurs'' are exactly the same as
the parameters for which the reconstruction problem is not solvable.
More formally, for determining if the glassy phase occurs for random
$(\deg+1)$-regular graphs and Gibbs measures with some parameters, one needs
to check if the reconstruction problem for the $\deg$-ary tree and
associated parameters is solvable or not.

Furthermore, it is claimed in~\cite{MezardMontanari:06} that the
reconstruction problem determines the type of glassy phase as
follows.  {  Mezard and Montanari predict that one-step
replica symmetry breaking occurs exactly when
when the Kesten-Stigum bound is not equal to the
reconstruction bound; otherwise full replica symmetry breaking
{ occurs}.}

{ Thus our results proved here, in conjunction with the theoretical physics
predictions of~\cite{MezardMontanari:06},
suggest the existence of two types of glassy phases
for spin systems on random graphs.}  It is an interesting challenge to state
these predictions in a rigorous mathematical way and to prove or disprove
them.

\section{Preliminaries and General Result}

For convenience, we sometimes write the channel
\begin{displaymath}
M
=
\left(\begin{array}{cc}
1-\eps^+ & \eps^+\\
1-\eps^- & \eps^-
\end{array}\right).
\end{displaymath}
Note first that the
stationary distribution $\pi = (\pi_+,\pi_-)$ of $M$ is given by
\begin{equation*}
\pi_+ = \frac{1 - \eps^-}{1-\theta} = \frac{1}{2} - \frac{\delta}{2(1-\theta)}, \qquad
\pi_- = \frac{\eps^+}{1-\theta} = \frac{1}{2} + \frac{\delta}{2(1-\theta)}.
\end{equation*}
In particular, this expression implies that the stationary distribution depends
only on the ratio $\delta/(1-\theta)$. Or put differently, each two of
the parameters $\pi_+,\delta$ and $\theta$ determine the third one
uniquely.
Note also that
\begin{eqnarray*}
\theta = \eps^- - \eps^+, \qquad \pi_- - \pi_+ = \frac{\delta}{1 - \theta}.
\end{eqnarray*}
Without loss of generality, we assume throughout that $\pi_- \geq \pi_+$ or equivalently
that $\delta \geq 0$. (Note that $\delta$ can be made negative by inverting the role of $+$ and $-$.)
Below, we will use the notation
\begin{eqnarray*}
\pimp \equiv \pi_- \pi_+^{-1},\qquad \deltat \equiv \pi_{-/+} - 1.
\end{eqnarray*}

\subsection{General Trees}\label{sec:generaltrees}

In this section, we state our Theorem in a more general setting. Namely, we consider
general rooted trees where different edges are equipped with different
transition matrices---all having the same stationary distribution
$\pi = (\pi_+,\pi_-)$. In other words, we consider
a general infinite rooted tree $\tcal=(V,E)$
equipped with a function $\theta : E \to [-1,1]$ such
that the edge $e$ of the tree is equipped with the matrix $M^e$ with
$\theta(M^e) = \theta(e)$ and the stationary distribution of
$M^e$ is $(\pi_+,\pi_-)$.

In this general setting the notion of degree is extended to the notion
of {\em branching number}.
In \cite{Furstenberg:70}, Furstenberg introduced the Hausdorff
dimension of a tree. Later, Lyons \cite{Lyons:89,Lyons:90}
showed that many probabilistic properties of the tree are determined by this number
which he named the branching number.
 For our purposes it is best to define the branching number via
cutsets.
\begin{definition}[Cutsets]
A {\em  cutset}
$S$ for a tree $\tcal$ rooted at $\rho$, is a finite set of
vertices separating $\rho$ from $\infty$. In other words, a finite set
$S$ is a cutset if every infinite self avoiding path from $\rho$
intersects $S$. An {\em antichain { or minimal cutset}} is a cutset that does not have
any proper subset which is also a cutset.
\end{definition}
%
%
\begin{definition} [Branching Number]\label{def:braT}
Consider a rooted tree $\tcal=(V,E,\rho)$ equipped with an edge function
$\theta : E \to [-1,1]$. For each vertex $v \in V$ we define
\[
\eta(x) = \prod_{e \in \path(\rho,x)} \theta^2(e),
\]
where $\path(\rho, x)$ is the set of edges on the unique path
between $\rho$ and $x$ in $\tcal$. The branching number $\br(\tcal, \theta)$ of $(\tcal,\theta)$ is
defined as
\[
\br(\tcal, \theta) = \inf \left\{ \lam > 0 :
         \inf_{\cutsets\ S} \sum_{x \in S} \eta(x) \lam^{-|x|} = 0 \right\}.
\]
\end{definition}

In our main result we show
\begin{theorem} [Reconstructibility on General Trees]\label{thm:d}
Let $0 \leq \theta_0 < 1$. Then there exists $\delta_0 > 0$ such that, for
all distributions $\pi= (\pi_+,\pi_-)$ with $\max\{|\delta(\pi,\theta_0)|,|\delta(\pi,-\theta_0)|\} <
\delta_0$ and for all trees $(\tcal,\theta)$ with $\sup_e |\theta(e)| \leq
\theta_0$ and $\br(\tcal,\theta) \leq 1$,
the reconstruction problem is not solvable.
\end{theorem}
It is easy to see that the conditions of Theorem~\ref{thm:d} hold for
$T_{\deg}$ if $\theta(e) = \theta$ for all $e$ and
$\deg \theta^2 \leq 1$.

\subsection{Magnetization}\label{sec:magnetization}

Let $T$ be a finite tree rooted at $x$ with edge function $\theta$.
Let $\sigma$ be the leaf states generated by
the Markov chain on $(T, \theta)$ with
stationary distribution $(\pi_+,\pi_-)$.
We denote by $\prob^+_T, \expec^+_T$
(resp. $\prob^-_T, \expec^-_T$, and $\prob_T, \expec_T$)
the probability/expectation operators with respect to the measure on the leaves of $T$
obtained by conditioning the root to be $+$ (resp. $-$, and stationary).
With a slight abuse of notation, we also write $\prob_T[+\,|\,\sigma]$
for the probability that the state at the the root of $T$ is $+$ given state
$\sigma$ at the leaves. The main random variable we consider is the
\emph{weighted magnetization of the root}
\begin{equation*}
X = \pi_-^{-1}\left[\pi_-\prob_T[+\,|\,\sigma] - \pi_+ \prob_T[-\,|\,\sigma]\right].
\end{equation*}
Note that the weights are chosen to guarantee
\begin{equation*}
\expec_T[X] = \pi_-^{-1}\left[\pi_- \pi_+ - \pi_+\pi_- \right]
= 0,
\end{equation*}
while the factor $\pi_-^{-1}$ is such that $|X|\leq 1$ with
probability $1$.

{ Note that for any random variable depending only on the
leaf states,
 $f = f(\sigma)$, we have
 $\pi_+ \expec^+_T[f] + \pi_- \expec^-_T[f] = \expec_T[f]$,
 so that in particular}
\begin{equation*}
\pi_+ \expec^+_T[X] + \pi_- \expec^-_T[X] = \expec_T[X] = 0,
\qquad
\pi_+ \expec^+_T[X^2] + \pi_- \expec^-_T[X^2] = \expec_T[X^2].
\end{equation*}
We define the following analogues of the Edwards-Anderson
order parameter for spin glasses on trees rooted at $x$
\begin{equation*}
\mom{x} = \expec_T[X^2],\qquad
\momp{x} = \expec^+_T[X^2],\qquad
\momm{x} = \expec^-_T[X^2].
\end{equation*}

\medskip

Now suppose $\tcal$ is an infinite tree rooted at $\rho$ with edge function $\theta$.
Let $T_n = (V_n,E_n,x_n)$, where
$V_n = \{u \in V: d(u,\rho) \leq n\}$,{
$E_n = \{e \in E:  d(e,\rho) \leq n\}$, and $x_n$ is identified with $\rho$.
It is not hard} to see that non-reconstructibility on $(T,\theta)$ is equivalent
in our notation to
\begin{eqnarray*}
\limsup_{n\to \infty} \mom{x}_n = 0.
\end{eqnarray*}
(Note that the total variation distance is
monotone in the cutsets. Therefore the limit goes to 0 with the
levels if and only if there exists a sequence of cutsets for which it goes to 0.)

\subsection{Expectations}\label{sec:expectations}

\begin{figure}
\centering
\input{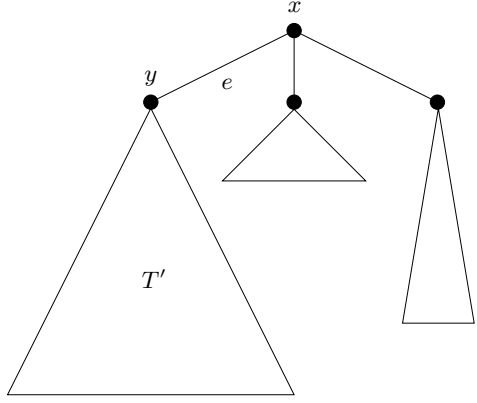}\caption{A finite tree $T$.}
\label{fig:fig1}
\end{figure}

Fix a stationary distribution $\pi = (\pi_+, \pi_-)$.
Let $T=(V,E)$ be a finite tree rooted at $x$ with edge function
$\{\theta(f), f\in E\}$
and weighted magnetization at the root $X$.
Let $y$ be a child of $x$ and $T'$ be the subtree of $T$ rooted at $y$.
Let $Y$ be the weighted magnetization at the root of $T'$. See Figure~\ref{fig:fig1}.
Denote by $\sigma$ the leaf states of $T$ and let $\sigma'$ be the restriction
of $\sigma$ to the leaves of $T'$. Assume the channel on $e = (x,y)$
is given by
\begin{displaymath}
M^e
=
\left(\begin{array}{cc}
1-\eps^+ & \eps^+\\
1-\eps^- & \eps^-
\end{array}\right)
=
\frac{1}{2}\left[
\left(\begin{matrix}
1+\theta & 1-\theta\\
1-\theta & 1+\theta
\end{matrix}\right)
+
\delta\left(\begin{matrix}
-1 & 1\\
-1 & 1
\end{matrix}\right)
\right].
\end{displaymath}
We collect in the next lemmas a number of useful identities.
\begin{lemma}[Radon-Nikodym Derivative]\label{lem:radon}
The following hold:
\begin{eqnarray*}
\frac{\diff \prob^+_T}{\diff \prob_T} = 1 + \pimp X,
\qquad
\frac{\diff \prob^-_T}{\diff \prob_T} = 1 - X,
\end{eqnarray*}
\begin{eqnarray*}
\expec^+_T[X] = \pimp \expec_T[X^2], \qquad \expec^-_T[X] = - \expec_T[X^2].
\end{eqnarray*}
\end{lemma}
\begin{proof}
Note that
\begin{eqnarray*}
X
= \pi_-^{-1} \left[\pi_- \prob_T[+\,|\,\sigma] - \pi_+\prob_T[-\,|\,\sigma]\right]
= \pi_-^{-1} \left[\prob_T[+\,|\,\sigma] - \pi_+\right]
= \pimp^{-1} \left[\frac{\prob_T[+\,|\,\sigma]}{\pi_+} - 1\right],
\end{eqnarray*}
so that
\begin{equation*}
\frac{\diff \prob^+_T}{\diff \prob_T} = \frac{\prob_T[+\,|\,\sigma]}{\pi_+} = 1 + \pimp X.
\end{equation*}
Likewise,
\begin{equation*}
\frac{\diff \prob^-_T}{\diff \prob_T} = \frac{\prob_T[-\,|\,\sigma]}{\pi_-} = 1 - X.
\end{equation*}

\medskip

Then, it follows that
\begin{equation*}
\expec^+_T[X]
=
\expec_T\left[X\left(1 + \pimp X\right)\right]
= \pimp \expec_T[X^2],
\end{equation*}
and similarly for $\expec^-_T[X]$.
\end{proof}
\begin{lemma}[Child Magnetization]\label{lem:child}
We have,
\begin{eqnarray*}
\expec^+_T[Y] = \theta \expec^+_{T'}[Y],\qquad \expec^-_T[Y] = \theta \expec^-_{T'}[Y],
\end{eqnarray*}
and
\begin{eqnarray*}
\expec^+_T[Y^2] = (1-\theta)\expec_{T'}[Y^2] + \theta \expec^+_{T'}[Y^2],
\qquad \expec^-_{T}[Y^2] = (1-\theta)\expec_{T'}[Y^2] + \theta \expec^-_{T'}[Y^2].
\end{eqnarray*}
\end{lemma}
\begin{proof}
By the Markov property, we have
\begin{eqnarray*}
\expec^+_T[Y]
&=&
(1-\eps^+)\expec^+_{T'}[Y] + \eps^+\expec^-_{T'}[Y]
=
\left[(1-\eps^+) - \eps^+\frac{\pi_+}{\pi_-}\right] \expec^+_{T'}[Y]
=
\left[(1-\eps^+) - (1-\eps^-)\right] \expec^+_{T'}[Y]\\
&=&
\theta \expec^+_{T'}[Y],
\end{eqnarray*}
and similarly for $\expec^-_T[Y]$.

\medskip

Also,
\begin{eqnarray*}
\expec^+_{T}[Y^2]
&=&
(1-\eps^+)\expec^+_{T'}[Y^2] + \eps^+\expec^-_{T'}[Y^2]
=
(1-\eps^+)\expec^+_{T'}[Y^2] + \frac{\eps^+}{\pi_-}(\expec_{T'}[Y^2] - \pi_+\expec^+_{T'}[Y^2])\\
&=&
\theta \expec^+_{T'}[Y^2] + (1-\theta)\expec_{T'}[Y^2],
\end{eqnarray*}
where we have used the calculation above. A similar expression holds for
$\expec^-_{T}[Y^2]$.
\end{proof}

\section{Tree Operations}\label{sec:recursion}

\begin{figure}
\centering
\input{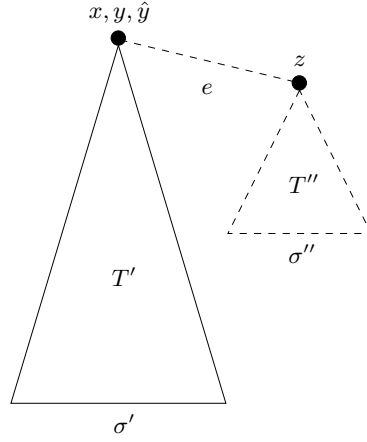}\caption{Tree $T$ after the \emph{Add-Merge} of $T'$ and $T''$. The dashed subtree is $\that$.}
\label{fig:fig2}
\end{figure}

To derive moment recursions, the basic graph operation we perform is the
following \emph{Add-Merge} operation.
Fix a stationary distribution $\pi = (\pi_+, \pi_-)$.
Let $T'$ (resp. $T''$) be a finite tree rooted at $y$ (resp. $z$) with
edge function $\theta'$ (resp. $\theta''$), leaf state
$\sigma'$ (resp. $\sigma''$), and weighted magnetization at the root $Y$
(resp. $Z$). Now add an edge $e = (\yhat,z)$ with edge value $\theta(e) = \theta$
to $T''$ to obtain a new tree $\that$. Then merge $\that$ with $T'$
by identifying $y = \yhat$ to obtain a new tree $T$. To avoid ambiguities, we denote by $x$ the root of $T$ and
$X$ the magnetization of the root of $T$ (where we identify the edge function on $T$ with
those on $T'$, $T''$, and $e$).  We let
$\sigma = (\sigma', \sigma'')$ be the leaf state of $T$. See Figure~\ref{fig:fig2}.
Let also $\yyhat$ be the magnetization of the root on $\that$.
Assume
\begin{displaymath}
M^e
=
\left(\begin{array}{cc}
1-\eps^+ & \eps^+\\
1-\eps^- & \eps^-
\end{array}\right).
\end{displaymath}
We first analyze the effect of adding an edge and merging subtrees on the magnetization variable.
\begin{lemma}[Adding an Edge]\label{lem:adding}
With the notation above, we have
\begin{equation*}
\yyhat = \theta Z.
\end{equation*}
\end{lemma}
\begin{proof}
Note that by Bayes' rule, the Markov property, and Lemma~\ref{lem:radon},
\begin{eqnarray*}
\yyhat
&=& \pi_+ \sum_{\gamma = +,-} \gamma\, \frac{\prob_{\that}[\gamma\,|\,\sigma'']}{\pi^\gamma}
= \pi_+\sum_{\gamma = +,-} \gamma\, \frac{\prob_{\that}[\sigma''\,|\,\gamma]}{\prob_{\that}[\sigma'']}\\
&=& \pi_+ \frac{\prob_{T''}[\sigma'']}{\prob_{\that}[\sigma'']}
\sum_{\gamma = +,-}\gamma\,
\left[(1-\eps^\gamma)\frac{\prob_{T''}[\sigma''\,|\,+]}{\prob_{T''}[\sigma'']}
+ \eps^\gamma\frac{\prob_{T''}[\sigma'\,|\,-]}{\prob_{T''}[\sigma'']}\right]\\
&=& \pi_+
\sum_{\gamma = +,-}\gamma\,
\left[(1-\eps^\gamma)\left(1 + \pimp Z\right)
+ \eps^\gamma \left(1 - Z\right)\right],
\end{eqnarray*}
where we have used $\prob_{\that}[\sigma''] = \prob_{T''}[\sigma'']$.
We now compute the expression in square brackets. We have
\begin{equation*}
(1-\eps^\gamma)\left(1 + \pimp Z\right)
+ \eps^\gamma \left(1 - Z\right)
= 1 + \pi_- Z\left[\frac{1 - \eps^\gamma}{\pi_+} - \frac{\eps^\gamma}{\pi_-}\right].
\end{equation*}
For $\gamma = +$, we get
\begin{equation*}
\frac{1 - \eps^+}{\pi_+} - \frac{\eps^+}{\pi_-}
= (1-\theta)\left[\frac{1-\eps^+}{1-\eps^-} - 1\right]
= (1-\theta)\left[\frac{\eps^--\eps^+}{1-\eps^-}\right]
= \frac{\theta}{\pi_+}.
\end{equation*}
A similar calculation for the $-$ case gives for $\gamma=+,-$
\begin{equation*}
(1-\eps^\gamma)\left(1 + \pimp Z\right)
+ \eps^\gamma \left(1 - Z\right)
= 1 + \gamma\theta\pi_-\pi_\gamma^{-1}Z.
\end{equation*}
Plugging above gives $\yyhat = \theta Z$.
\end{proof}

\begin{lemma}[Merging Subtrees]\label{lem:merging}
With the notation above, we have
\begin{equation*}
X = \frac{Y + \yyhat + \deltat  Y \yyhat}
{1 + \pimp Y \yyhat}.
\end{equation*}
The same expression holds for a general $\that$.
\end{lemma}
\begin{proof}
By Bayes' rule, the Markov property, and Lemma~\ref{lem:radon}, we have
\begin{eqnarray*}
X
&=& \pi_+ \sum_{\gamma = +,-} \gamma\, \frac{\prob_T[\gamma\,|\,\sigma]}{\pi^\gamma}
= \pi_+\sum_{\gamma = +,-} \gamma\, \frac{\prob_T[\sigma\,|\,\gamma]}{\prob_T[\sigma]}
= \pi_+ \frac{\prob_{T'}[\sigma']\prob_{\that}[\sigma'']}{\prob_T[\sigma]}
\sum_{\gamma = +,-}\gamma\, \frac{\prob_{T'}[\sigma'\,|\,\gamma]}{\prob_{T'}[\sigma']}
\frac{\prob_{\that}[\sigma''\,|\,\gamma]}{\prob_{\that}[\sigma'']}\\
&=& \pi_+ \frac{\prob_{T'}[\sigma']\prob_{\that}[\sigma'']}{\prob_T[\sigma]}
\sum_{\gamma = +,-}\gamma\, \left[1 + \gamma \pi_{-} \pi_\gamma^{-1} (Y + \yyhat)
+ (\pi_{-} \pi_\gamma^{-1})^2 Y \yyhat \right].
\end{eqnarray*}
Similarly, we have
\begin{eqnarray*}
\frac{\prob_T[\sigma]}{\prob_{T'}[\sigma']\prob_{\that}[\sigma'']}
&=& \frac{1}{\prob_{T'}[\sigma']\prob_{\that}[\sigma'']} \sum_{\gamma = +,-} \pi^\gamma\, \prob_T[\sigma\,|\,\gamma]
= \sum_{\gamma = +,-}\pi_\gamma\,
\left[1 + \gamma \pi_{-} \pi_\gamma^{-1} (Y + \yyhat)
+ (\pi_{-} \pi_\gamma^{-1})^2 Y \yyhat \right].
\end{eqnarray*}
Note that
\begin{equation*}
\sum_{\gamma = +,-}\gamma\, \left[1 + \gamma \pi_{-} \pi_\gamma^{-1} (Y + \yyhat)
+ (\pi_{-} \pi_\gamma^{-1})^2 Y \yyhat \right]
=
\pi_+^{-1} (Y + \yyhat)
+ \pi_+^{-2}(\pi_- - \pi_+)Y \yyhat,
\end{equation*}
where we have used
\begin{eqnarray*}
\pi_-^2 - \pi_+^2 = (\pi_- - \pi_+)(\pi_- + \pi_+) = \pi_- - \pi_+.
\end{eqnarray*}
Similarly,
\begin{equation*}
\sum_{\gamma = +,-}\pi_\gamma\,
\left[1 + \gamma \pi_{-} \pi_\gamma^{-1} (Y + \yyhat)
+ (\pi_{-} \pi_\gamma^{-1})^2 Y \yyhat \right]
=
1
+ \pi_- \pi_+^{-1} Y \yyhat.
\end{equation*}
The result follows.
\end{proof}

\section{Symmetric Channels On Regular Trees}

As a warm-up, we start by analyzing the binary symmetric channel on the infinite $\deg$-ary tree.
Our proof is arguably the simplest proof to date of this result. The same proof structure will be
used in the general case.
\begin{theorem}[Symmetric Channel. See \cite{BlRuZa:95, EvKePeSc:00,Ioffe:96a, MosselPeres:03, JansonMossel:04, BeKeMoPe:05,MaSiWe:04}.]\label{thm:symmetric}
Let $M$ be a transition matrix with $\delta = 0$ and
$\deg\theta^2 \leq 1$. Let $\tcal$ be the infinite $\deg$-ary tree. Then, the
reconstruction problem on $(\tcal,M)$ is not solvable.
\end{theorem}
\begin{proof}
Consider again the setup of Section~\ref{sec:recursion}.
Note first that, by Lemmas~\ref{lem:radon},~\ref{lem:child} and~\ref{lem:adding}, we have
\begin{equation}\label{eq:first}
\expec_{\that}[\yyhat^2] = \expec^+_{\that}[\yyhat]
= \theta \expec^+_{\that}[Z]
= \theta^2 \expec^+_{T''}[Z] = \theta^2 \expec_{T''}[Z^2],
\end{equation}
where we have used the fact that $\pimp = 1$ when $\delta = 0$ (although note that
it is not needed).
In other words, adding an edge to the root of a tree and re-rooting at the new vertex
has the effect of multiplying the second moment of the magnetization by $\theta^2$.
Now consider the \emph{Add-Merge} operation
defined in Section~\ref{sec:recursion}.
Using the expansion
\begin{equation}\label{eq:expansion}
\frac{1}{1+r} = 1 - r + \frac{r^2}{1+r},
\end{equation}
the inequality $|X|\leq 1$,
and Lemma~\ref{lem:merging},
we get
\begin{eqnarray}\label{eq:ineq1}
X
= Y + \yyhat
- Y \yyhat (Y + \yyhat) + Y^2 \yyhat^2 X
\leq Y + \yyhat
- Y \yyhat (Y + \yyhat) + Y^2 \yyhat^2.
\end{eqnarray}
Note that from Lemmas~\ref{lem:radon} and~\ref{lem:child}, we have
\begin{eqnarray*}
\expec^+_{T}[X] = { \mom{x}},\quad
\expec^+_{T}[Y] = \expec^+_{T}[Y^2] = \mom{y},\quad
\expec^+_{T}[\yyhat] =
\expec^+_{T}[\yyhat^2] = \theta^2 \mom{z},
\end{eqnarray*}
where we have used that $\momp{y} = \momm{y} = \mom{y}$ and
$\momp{z} = \momm{z} = \mom{z}$ by symmetry.
Taking $\expec^+_{T}$ on both sides of (\ref{eq:ineq1}), we get
\begin{eqnarray*}
\mom{x} &\leq& \mom{y} + \theta^2 \mom{z} - \theta^2 \mom{y} \mom{z} - \theta^2 \mom{y} \mom{z} + \theta^2 \mom{y} \mom{z}
= \mom{y} + \theta^2 \mom{z} - \theta^2 \mom{y}\mom{z}.
\end{eqnarray*}

\medskip

Now, let $T_n = (V_n, E_n, x_n)$ be as in Section~\ref{sec:magnetization}.
Repeating the \emph{Add-Merge} operation $(\deg - 1)$ times, we finally have by induction
\begin{eqnarray*}
\mom{x}_n \leq \deg \theta^2 \mom{x}_{n-1} - (\deg-1) \theta^4 \mom{x}_{n-1}^2.
\end{eqnarray*}
Indeed, note that for $0< a < \deg$,
\begin{eqnarray*}
(a \theta^2 \mom{x}_{n-1} - (a-1) \theta^4 \mom{x}_{n-1}^2) + \theta^2 \mom{x}_{n-1}
- \theta^2 (a \theta^2 \mom{x}_{n-1} - (a-1) \theta^4 \mom{x}_{n-1}^2) \mom{x}_{n-1}
\leq (a+1) \theta^2 \mom{x}_{n-1} - a \theta^4 \mom{x}_{n-1}^2,
\end{eqnarray*}
and the first step of the induction is given by~(\ref{eq:first}).
This concludes the proof.
\end{proof}

\section{Roughly Symmetric Channels on General Trees}

We now tackle the general case. We start by analyzing the \emph{Add-Merge}
operation.
\begin{proposition} [Basic Inequality]\label{prop:basic}
Consider the setup of Section~\ref{sec:recursion}. Assume $|\theta| < 1$.
Then, there is a $\delta_0(|\theta|) > 0$ depending
only on $|\theta|$ such that
\begin{eqnarray*}
\mom{x} \leq \mom{y} + \theta^2 \mom{z},
\end{eqnarray*}
whenever $\delta$ (on $e$) is  less than $\delta_0(|\theta|)$.
\end{proposition}
\begin{proof}
The proof is similar to that in the symmetric case.
By expansion (\ref{eq:expansion}), inequality $|X|\leq 1$,
and Lemma~\ref{lem:merging}, we have
\begin{eqnarray}\label{eq:ineq2}
X
\leq Y + \yyhat + \Delta Y \yyhat
- \pimp Y \yyhat (Y + \yyhat + \Delta Y \yyhat)
+ \pimp^2 Y^2 \yyhat^2.
\end{eqnarray}
Let $\rho' = (\mom{y})^{-1} \momp{y}$ and $\rho'' = (\mom{z})^{-1} \momp{z}$. Then,
by Lemmas~\ref{lem:radon} and~\ref{lem:child}, we have
\begin{eqnarray*}
&\expec^+_{T}[X] = \pimp \mom{x},\quad
\expec^+_{T}[Y] = \pimp \mom{y},\quad
\expec^+_{T}[Y^2] = \mom{y}\rho',\quad\\
&\expec^+_{T}[\yyhat] = \pimp \theta^2 \mom{z},\quad
\expec^+_{T}[\yyhat^2] = \theta^2 \mom{z} [(1-\theta) + \theta \rho''].
\end{eqnarray*}
Taking $\pimp^{-1} \expec^+_{T}$ on both sides of (\ref{eq:ineq2}), we get
\begin{eqnarray*}
\mom{x} &\leq&
\mom{y} + \theta^2 \mom{z}
+ \Delta \pimp \theta^2 \mom{y}\mom{z}\\
&& \quad - \pimp \theta^2 \mom{y}\mom{z} \rho'
- \pimp \theta^2 \mom{y}\mom{z} [(1-\theta) + \theta \rho'']
- \Delta \theta^2 \mom{y}\mom{z} \rho' [(1-\theta) + \theta \rho'']\\
&& \quad\quad + \pimp \theta^2 \mom{y}\mom{z} \rho' [(1-\theta) + \theta \rho'']\\
&\leq& \mom{y} + \theta^2 \mom{z} - \pimp \theta^2 \mom{y}\mom{z}[\acal - \Delta \bcal],
\end{eqnarray*}
where
\begin{eqnarray*}
\acal = \rho' + (1 - \rho')[(1-\theta) + \theta \rho''],
\end{eqnarray*}
and
\begin{eqnarray*}
\bcal = 1 - \pimp^{-1} \rho' [(1-\theta) + \theta \rho''].
\end{eqnarray*}
Note that $[(1-\theta) + \theta \rho''] \geq 0$ by Lemma~\ref{lem:child}.
So $\bcal \leq 1$ and it suffices to have $\acal \geq \Delta$.
Note also that $\acal$ is multilinear in $(\rho', \rho'')$.
Therefore, to minimize $\acal$, we only need to consider extreme cases in $(\rho',\rho'')$.
By $\pi_+ y^+ + \pi_- y^- = y$ it follows that $0\leq \rho'\leq \pi_+^{-1}$. The same
holds for $\rho''$. At $\rho' = 0$, we have
\begin{displaymath}
\acal = 1 - \theta[1 - \rho'']
\geq
\left\{
\begin{array}{ll}
1 - \theta, & \mathrm{if}\ \theta \geq 0,\\
1 - \pimp |\theta|, & \mathrm{if}\ \theta \leq 0,
\end{array}
\right.
\end{displaymath}
where we have used
\begin{eqnarray*}
1 - \pi_+^{-1} = - \pimp.
\end{eqnarray*}
At $\rho' = \pi_+^{-1}$, we have
\begin{displaymath}
\acal = \pi_+^{-1} + (1 - \pi_+^{-1})[1 - \theta[1 - \rho'']]
= 1 + \theta \pimp [1 - \rho'']
\geq
\left\{
\begin{array}{ll}
1 - \pimp^2 \theta, & \mathrm{if}\ \theta \geq 0,\\
1 - \pimp |\theta|, & \mathrm{if}\ \theta \leq 0.
\end{array}
\right.
\end{displaymath}
Since $\pimp \geq 1$ by assumption, it follows that
\begin{eqnarray*}
\acal \geq 1 - \pimp^2 |\theta|.
\end{eqnarray*}
At $\delta = 0$, this bound is strictly positive and moreover $\Delta = 0$. Therefore,
by continuity in $\delta$ of $\Delta$ and the bound above, the result
follows.
\end{proof}

\begin{proposition}[Induction Step]\label{prop:induction}
Let $T$ be a finite tree rooted at $x$ with edge function $\theta$.
Let $w_1,\ldots,w_\alpha$ be the children of $x$ in $T$ and denote
by $e_a$ the edge connecting $x$ to $w_a$. Let
$\theta_0 = \max\{|\theta(e_1)|,\ldots,|\theta(e_\alpha)|\}$ and
assume that on each edge $e_a$, $\delta \leq \delta_0(\theta_0)$,
where $\delta_0$ is defined in Proposition~\ref{prop:basic}. Then
\begin{equation*}
\mom{x} \leq \sum_{a = 1}^\alpha \theta(e_a)^2 \mom{w}_a.
\end{equation*}
\end{proposition}
\begin{proof}
As noted in the proof of Theorem~\ref{thm:symmetric},
adding an edge $e$ to the root of a tree and re-rooting at the new vertex
has the effect of multiplying the second moment of the magnetization by $\theta^2(e)$.
The result follows by applying Proposition~\ref{prop:basic} $(\alpha - 1)$ times.
\end{proof}

\noindent{\bf Proof of Theorem~\ref{thm:d}}:
It suffices to show that for all $\eps > 0$ there is an $N$ large enough so that
$\mom{x}_n \leq \eps$, $\forall n \geq N$. Fix $\eps > 0$. By definition of the branching
number, there exists a cutset $S$ of $\tcal$ such that
\begin{eqnarray*}
 \sum_{u \in S} \eta(u) \leq \eps.
 \end{eqnarray*}
Assume w.l.o.g. that $S$ is actually an antichain and let $N$ be such that $S$ is in $T_N$.
It is enough to show that
\begin{eqnarray*}
\mom{x}_n \leq \sum_{u \in S} \eta(u),\qquad \forall\ n \geq N.
\end{eqnarray*}
Fix $n \geq N$. Applying Proposition~\ref{prop:induction} repeatedly from the root of $T_n$ down
to $S$, it is clear that
\begin{eqnarray*}
\mom{x}_n \leq \sum_{u \in S} \eta(u) \expec_{T_n(u)}[U^2]  \leq \sum_{u \in S} \eta(u),
\end{eqnarray*}
where $T_n(u)$ is the subtree of $T_n$ rooted at $u$ and $U$ is the
magnetization at $u$ on $T_n(u)$ (with $|U| \leq 1$). This concludes the proof.
$\blacksquare$

\clearpage

\bibliographystyle{plain}
\bibliography{all,myshort}

\begin{thebibliography}{10}

\bibitem{BeKeMoPe:05}
N.~Berger, C.~Kenyon, E.~Mossel, and Y.~Peres.
\newblock Glauber dynamics on trees and hyperbolic graphs.
\newblock {\em Probab. Theory Related Fields}, 131(3):311--340, 2005.
\newblock Extended abstract by Kenyon, Mossel and Peres appeared in proceedings
  of 42nd IEEE Symposium on Foundations of Computer Science (FOCS) 2001,
  568--578.

\bibitem{BlRuZa:95}
P.~M. Bleher, J.~Ruiz, and V.~A. Zagrebnov.
\newblock On the purity of the limiting {G}ibbs state for the {I}sing model on
  the {B}ethe lattice.
\newblock {\em J. Statist. Phys.}, 79(1-2):473--482, 1995.

\bibitem{CCST:86}
J.~T. Chayes, L.~Chayes, James~P. Sethna, and D.~J. Thouless.
\newblock A mean field spin glass with short-range interactions.
\newblock {\em Comm. Math. Phys.}, 106(1):41--89, 1986.

\bibitem{DaMoRo:06}
C.~Daskalakis, E.~Mossel, and S.~Roch.
\newblock {Optimal Phylogenetic Reconstruction}.
\newblock Availible on the Arxiv at math.PR/0509575, Extended abstract to
  appear at Proceedings of STOC 2006, 2006.

\bibitem{EvKePeSc:00}
W.~S. Evans, C.~Kenyon, Yuval Y.~Peres, and L.~J. Schulman.
\newblock Broadcasting on trees and the {I}sing model.
\newblock {\em Ann. Appl. Probab.}, 10(2):410--433, 2000.

\bibitem{Felsenstein:04}
J.~Felsenstein.
\newblock {\em Inferring Phylogenies}.
\newblock Sinauer, New York, New York, 2004.

\bibitem{Furstenberg:70}
H.~Furstenberg.
\newblock Intersections of {C}antor sets and transversality of semigroups.
\newblock In {\em Problems in analysis (Sympos. Salomon Bochner, Princeton
  Univ., Princeton, N.J., 1969)}, pages 41--59. Princeton Univ. Press,
  Princeton, N.J., 1970.

\bibitem{Georgii:88}
H.~O. Georgii.
\newblock {\em Gibbs measures and phase transitions}, volume~9 of {\em de
  Gruyter Studies in Mathematics}.
\newblock Walter de Gruyter \& Co., Berlin, 1988.

\bibitem{Higuchi:77}
Y.~Higuchi.
\newblock Remarks on the limiting {G}ibbs states on a {$(d+1)$}-tree.
\newblock {\em Publ. Res. Inst. Math. Sci.}, 13(2):335--348, 1977.

\bibitem{Ioffe:96a}
D.~Ioffe.
\newblock On the extremality of the disordered state for the {I}sing model on
  the {B}ethe lattice.
\newblock {\em Lett. Math. Phys.}, 37(2):137--143, 1996.

\bibitem{JansonMossel:04}
S.~Janson and E.~Mossel.
\newblock Robust reconstruction on trees is determined by the second
  eigenvalue.
\newblock {\em Ann. Probab.}, 32:2630--2649, 2004.

\bibitem{KestenStigum:66}
H.~Kesten and B.~P. Stigum.
\newblock Additional limit theorems for indecomposable multidimensional
  {G}alton-{W}atson processes.
\newblock {\em Ann. Math. Statist.}, 37:1463--1481, 1966.

\bibitem{Lyons:89}
R.~Lyons.
\newblock The {I}sing model and percolation on trees and tree-like graphs.
\newblock {\em Comm. Math. Phys.}, 125(2):337--353, 1989.

\bibitem{Lyons:90}
R.~Lyons.
\newblock Random walks and percolation on trees.
\newblock {\em Ann. Probab.}, 18(3):931--958, 1990.

\bibitem{MaSiWe:04}
F.~Martinelli, Alistair A.~Sinclair, and D.~Weitz.
\newblock Glauber dynamics on trees: boundary conditions and mixing time.
\newblock {\em Comm. Math. Phys.}, 250(2):301--334, 2004.

\bibitem{MezardZecchina:02}
M.~Mezard and R.~Zecchina.
\newblock Random k-satisfiability: from an analytic solution to an effici ent
  algorithm.
\newblock {\em Phys. Rev. E}, 66, 2002.

\bibitem{MezardMontanari:06}
M.~M\'{e}zard~A. Montanari.
\newblock Reconstruction on trees and the spin glass transition, 2006.
\newblock Preprint.

\bibitem{Mossel:01}
E.~Mossel.
\newblock Reconstruction on trees: beating the second eigenvalue.
\newblock {\em Ann. Appl. Probab.}, 11(1):285--300, 2001.

\bibitem{Mossel:04}
E.~Mossel.
\newblock Survey: Information flow on trees.
\newblock In J.~Nestril and P.~Winkler, editors, {\em Graphs, Morphisms and
  Statistical Physics. DIMACS series in discrete mathematics and theoretical
  computer science}, pages 155--170. Amer. Math. Soc., 2004.

\bibitem{MosselPeres:03}
E.~Mossel and Y.~Peres.
\newblock Information flow on trees.
\newblock {\em Ann. Appl. Probab.}, 13(3):817--844, 2003.

\bibitem{MePaZe:03}
M.~M\'{e}zard~G. Parisi and R.~Zecchina.
\newblock Analytic and algorithmic solution of random satisfiability problems.
\newblock {\em Science}, 297, 812, 2002.
\newblock (Scienceexpress published on-line 27-June-2002; 10.1126/science.
  1073287).

\bibitem{PemantlePeres:06}
Robin Pemantle and Yuval Peres.
\newblock {The critical Ising model on trees, concave recursions and nonlinear
  capacity}.
\newblock Available at: arXiv:math.PR/0503137.

\bibitem{Spitzer:75}
F.~Spitzer.
\newblock Markov random fields on an infinite tree.
\newblock {\em Ann. Probability}, 3(3):387--398, 1975.

\end{thebibliography}

\clearpage

\appendix
{
\section{Lower bound on $\delta_0$}\label{sec:appendix}

\begin{lemma}[Bound on $\delta_0$]
Let $\delta_0$ be as in Propostion~\ref{prop:basic}.
Let $0 \leq \theta_0 < 1$. Then, $\delta_0(\theta_0)$ can be set
as large as $\bar\delta = (1 - \theta_0)\beta(\theta_0)$,
where $\beta(\theta_0)$ is the smallest root of
\begin{eqnarray*}
(1 - \theta_0) - (4 + 2 \theta_0) \beta + (3 - \theta_0) \beta^2 = 0.
\end{eqnarray*}
In particular, if $\theta_0 = 1/\sqrt{\deg}$ (as in the $\deg$-ary case), 
$\bar\delta \approx 0.016$ when $\deg = 2$ and 
$\bar\delta \approx 1/3$ when $\deg$ is large.
\end{lemma}
\begin{proof}
Let
\begin{equation*}
\phi = \frac{\delta}{1 - \theta_0}.
\end{equation*} 
Then (letting $|\theta|=\theta_0$)
\begin{equation*}
\pimp = \frac{1 + \phi}{1 - \phi}.
\end{equation*}
From the proof of Proposition~\ref{prop:basic}, we seek the largest value of $\delta \geq 0$ such that
\begin{equation*}
(1 - \pimp^2 \theta_0) - (\pimp - 1) \geq 0.
\end{equation*}
Multiplying by $(1 - \phi)^2$ and rearranging, we get
\begin{equation*}
2(1 - \phi)^2 - (1+\phi)(1-\phi) - \theta_0 (1+\phi)^2
= (1 - \theta_0) - (4 + 2\theta_0)\phi + (3 - \theta_0)\phi^2.
\end{equation*}
This expression is positive at $\phi=0$ and remains positive
until it reaches its smallest root in $\phi$.

\medskip

When $\theta_0 = 0$, the polynomial above reduces to
\begin{equation*}
1 - 4\phi + 3\phi^2 = (1 - 3\phi)(1 - \phi),
\end{equation*}
which has its smallest root at $1/3$.
The special case $\deg = 2$ in the statement of the lemma can be computed
numerically.
\end{proof}
}
\end{document}